\documentclass[11pt]{amsart}
\usepackage{amsmath,amssymb,amsfonts,url,mathptmx}
\numberwithin{equation}{section}

\newtheorem{thm}{Theorem}[section]

\newtheorem{rem}{Remark}[section]

\newcommand{\hdot}{^\text{\r{}}\hspace{-.33cm}H}

\begin{document}
\title[Toda system]{Energy concentration and a priori estimates for $B_2$ and $G_2$ types of Toda systems} \subjclass{35J60, 35J47}
\keywords{SU(n+1)-Toda system, holomorphic curves, asymptotic analysis, a priori estimate, classification theorem, topological degree, blowup solutions}

\author{Chang-shou Lin}
\address{Taida Institute for Mathematical Sciences\\
        Center for the Advanced Study in Theoretical Sciences\\
        National Taiwan University\\
         Taipei 106, Taiwan } \email{cslin@math.ntu.edu.tw}

\author{Lei Zhang}
\address{Department of Mathematics\\
        University of Florida\\
        358 Little Hall P.O.Box 118105\\
        Gainesville FL 32611-8105}
\email{leizhang@ufl.edu}

\date{\today}

\begin{abstract} For Toda systems with Cartan matrix either $B_2$ or $G_2$, we prove that
the local mass of blowup solutions at its blowup points converges to a finite set. Further more this finite set can be completely determined for $B_2$ Toda systems, while for $G_2$ systems we need one additional assumption.  As an application of the local mass classification we establish a priori estimates for corresponding
 Toda systems defined on Riemann surfaces.
\end{abstract}


\maketitle

\section{Introduction}

Let $(M,g)$ be a Riemann surface with area equal to $1$ and Gauss curvature equal to a constant $K_0$. We consider the following Liouville equation
\begin{equation}\label{liou-e}
\Delta_g u+2 (e^u-K_0)=4\pi\sum_{j=1}^N\alpha_j \delta_{p_j} \quad \mbox{ in }\quad M
\end{equation}
where $\Delta_g$ is the Laplace-Beltrami operator ($-\Delta_g\ge 0$), $\alpha_j>- 1$ and $\delta_{p_j}$ is the Dirac measure at $p_j\in M$. The geometric meaning of (\ref{liou-e}) is that for any solution $u$, the new conformal metric
$ds^2=e^ug$ has the constant Gauss curvature equal to $1$ outside the singular points $\{p_j\}$. Near each $p_j$, using a complex coordinate $z$ satisfying $z(p_j)=0$ we have $e^{u(z)}=O(1)|z|^{2\alpha_j}$ for $|z|$ near $0$, hence , the new metric $ds^2$ is degenerate at $p_j$ and is called a metric with conic singularity. Equation (\ref{liou-e}) and its general form, the so-called mean field equation, have been extensively studied for many decades. For example see \cite{bm,chenlin1,chenlin2,prajapat,tr-1,tr-2,yang2} and the reference therein. In particular, for the case $M$ being the standard sphere or a torus, the equation can be written as ( replacing $u$ in (\ref{liou-e}) by $u+\log 2$)
\begin{equation}\label{liou-2}
\Delta_g u+e^u=4\pi \sum_j \alpha_j \delta_{p_j}.
\end{equation}
Recently equation (\ref{liou-2}) is found to have deep connection with the classical Lame equation and also the Painleve VI equation. For example for the Painleve VI equation with certain parameters, some non-existence theorem of (\ref{liou-2}) plays a key part in the proof of the smoothness of the solutions with unitary monodromy group. The interested readers may read into \cite{chai} and \cite{CKL} for more in-depth discussions.

The natural generalization of (\ref{liou-e}) is the so-called Toda system
\begin{equation}\label{toda-g}
\Delta_g u_i+\sum_{j=1}^n k_{ij} e^{u_j}=4\pi \sum_{j=1}^n \alpha_{ij} \delta_{p_{ij}}
\end{equation}
where $k_{ij}$ is known as the Cartan group of some simple Lie algebra. For the case of $A_m$ type (see explanations below), the coefficient matrix $A=(a_{ij})$ is expressed by
\begin{equation}\label{cartan-m}
A=\left(\begin{array}{ccccc}
2 & -1 & 0 & ... & 0 \\
-1 & 2 & -1 & ... & 0 \\
\vdots & \vdots &  &  & \vdots\\
0 & ... & -1 & 2 & -1 \\
0 & ... &  0 & -1 & 2
\end{array}
\right ).
\end{equation}
It is well known that the general Toda system (\ref{toda-g}) is closely related to geometry \cite{bolton1,bolton2,doliwa,guest} and the gauge theory in many physics models. For example, to describe the physics of high critical temperature superconductivity, a model of relative Chern-Simons model was proposed and this model can be reduced to a $n\times n$ system with exponential nonlinearity if the gauge potential and the Higgs field are algebraically restricted. Then the Toda system (\ref{toda-g}) with (\ref{cartan-m}) is one of the limiting equations if the coupling constant tends to zero. For extensive discussions on the relationship between Toda system and its background in Physics we refer the readers to \cite{bennet,dunne3,dunne2,dunne4,ganoulis,kao1,lee1,lee2,mansfield,yang2} and the reference therein.

If the rank of the simple Lie Algebra is $2$, there are three types of corresponding Cartan matrices of rank 2:
$$A_2=\left(\begin{array}{cc}
2 & -1 \\
-1 & 2
\end{array}
\right )
\quad B_2(=C_2)=\left(\begin{array}{cc}
 2 & -1 \\
 -2 & 2
 \end{array}
 \right )
 \quad G_2=\left(\begin{array}{cc}
 2 & -1 \\
 -3 & 2
 \end{array}
 \right ).
 $$

These rank $2$ matrices are the simplest examples of their more general forms. In general, there are four types of simple non-exceptional Lie Algebra: $A_m, B_m, C_m$ and $D_m$ whose Cartan subalgebra are
 $sl(m+1)$, $so(2m+1)$, $sp(m)$ and $so(2m)$, respectively.  Corresponding to each of the four types of Lie Algebra there is a Toda system. Solutions of Toda systems are closely related to holomorphic curves in projective spaces. In particular, from the classical Pl\"ucker formula we see that any holomorphic curve gives rise to a solution of the $A_m$ type Toda system and the branch points of these curves correspond to the singularities of the solutions. On the other hand, if we integrate the $A_m$ Toda system, a solution defines a holomorphic curve in $\mathbb C\mathbb P^n$ at least locally. The interested readers may see \cite{lin-wei-ye} for further discussions of this respect.

In \cite{chern} and \cite{uhlenbeck} it was noticed that the equation (\ref{toda-g}), like (\ref{liou-e}), is also an integrable system. The integrability has been further discussed in \cite{lin-wei-ye}. As we mentioned earlier, equation (\ref{toda-g}) has deep connections with algebraic geometry, modular forms and the Painleve VI equation. So it is natural for us to study (\ref{toda-g}) from analysis viewpoints as well as the perspectives from integrable systems. From the analytic side, the most important issue is to derive a degree counting formula for equation (\ref{toda-g}), a generalization of the previous works of Chen-Lin \cite{chenlin1,chenlin2} for (\ref{liou-e}) and Lin-Zhang \cite{linzhang1,linzhang2,lin-zhang-jfa} for general Liouville systems. However this generalization is very challenging because the bubbling phenomena are more complicated and the concentration has not yet been proven even for $SU(3)$ Toda system. So far in this direction the degree counting formula has only been proved for the simplest case of $SU(3)$ Toda system, see \cite{lin-wei-yang}.

In this paper we initiate the analytic program for the systems $B_2$ and $G_2$. The main purpose is to establish the a priori bound for non-critical parameters. Similar to (\ref{liou-e}) we consider the Toda system of the mean field type:

\begin{equation}\label{toda-n}
\Delta_g u_i+\sum_{j=1}^n a_{ij}\rho_j (\frac{h_je^{u_j}}{\int_M h_j e^{u_j}dV_g}-1)=0, i=1,..,n
\end{equation}
where $h_1,...,h_n$ are positive smooth functions on $M$, $\rho_1,...,\rho_n$ are positive constants, and the solutions are in the space
$$\hdot^{1,n}=\{v=(v_1,...,v_n);\quad \int_M v_i dV_g=0, v_i\in L^2(M),\,\, \nabla v_i\in L^2(M)\quad i=1,..,n\}. $$
For $B_2$ or $G_2$, the coefficient matrix $A$ is
$$A=\left(\begin{array}{cc}
2 & -1 \\
-2 & 2
\end{array}
\right )  \quad \mbox{ or }
\quad \left(\begin{array}{cc}
2 & -1 \\
-3 & 2
\end{array}
\right ), \mbox{ respectively}.
$$

One of the most important and intriguing issues of Toda system in general is the blowup phenomenon. A point $p$ is called a blowup point if, along a subsequence, a sequence of solutions $\{u^k=(u_1^k,...,u_n^k)\}$ satisfies
$$\max_{i}\max_{B(p,\delta)}u_i^k=\max_iu_i^k(p_k)\to \infty, \quad p_k\to p. $$
 Understanding the asymptotic behavior of blowup solutions near its blowup point is crucial for many important questions related to the Toda systems such as a priori estimates, degree counting formula, existence results and multiplicity results, etc. We say a sequence of blowup solutions $\{u^k\}$ possesses energy concentration if for some $i$, $\max u_i^k \to \infty$, $e^{u_i^k}$ tends to a Dirac measure as $k\to \infty$.  The purpose of this article is to study the blowup phenomenon and the energy concentration for (\ref{toda-n}) with $B_2$ or $G_2$ matrix.
For simplicity we assume

\begin{equation}\label{toda-b2g2}
\Delta u_i +\sum_{j=1}^2 a_{ij} h_j e^{u_j}=0, \quad \mbox{ in }\quad B(0,1),\quad i=1,2
\end{equation}
where $B(0,1)$ is the unit ball in $\mathbb R^2$ ( throughout the paper we use $B(p,r)$ to denote the ball centered at $p$ with radius $r$),
$A=(a_{ij})_{2\times 2}$ is a $B_2$ or $G_2$ matrix.

In more precise terms we let $u^k=(u_1^k,u_2^k)$ be a sequence of solutions to
\begin{equation}\label{toda-b2}
\Delta u_i^k +\sum_{j=1}^2 a_{ij} h_j^k e^{u_j^k}=0, \quad \mbox{ in }\quad B(0,1),\quad i=1,2
\end{equation}
where $h_1^k,h_2^k$ are two sequences of positive smooth functions with uniform bound:
\begin{equation}\label{ah}
\frac 1C\le h_i^k\le C, \quad \|h_i^k\|_{C^2(B(0,1))}\le C, \quad \mbox{in } B(0,1), \quad \forall i=1,2,
\end{equation}
the origin is the only blowup point for $u^k=(u_1^k,u_2^k)$, which has bounded energy and oscillation on $\partial B(0,1)$:
\begin{equation}\label{asm-u}
\left\{\begin{array}{ll}
\max_i\max_{B(0,1)}u_i^k\to \infty,\quad \max_i\max_{K\subset\subset B(0,1)\setminus \{0\}} u_i^k\le C(K),\\
\\
|u_i^k(x)-u_i^k(y)|\le C, \quad \forall x,y\in \partial B(0,1),\quad i=1,2.\\
\\
\int_{B(0,1)}h_i^k e^{u_i^k}\le C, \quad i=1,2.
\end{array}
\right.
\end{equation}
The main purpose of this article is to study the following quantity:
\begin{equation}\label{12sep2e3}
\sigma_i=\lim_{\delta\to 0}\lim_{k\to \infty} \frac{1}{2\pi} \int_{B_{\delta}} h_i^k e^{u_i^k},\quad i=1,2.
\end{equation}
Note that the inside limit: $\lim_{k\to \infty}$, which is to be taken first, is understood as taken along a subsequence of $u^k$ with the same notation. The second limit: $\lim_{\delta\to 0}$ indicates that we consider how the energy ( the integration of $h_i^ke^{u_i^k}$ ) concentrates at $0$.

Note that the oscillation finiteness assumption in (\ref{asm-u}) is natural and generally satisfied in most applications. The energy bound in (\ref{asm-u}) is also natural for systerm/equation defined in two dimensional spaces.

The first main result is:
\begin{thm}\label{cla-b2} Let $A$ be the $B_2$ matrix, the blowup solutions $u^k=(u_1^k,u_2^k)$ of (\ref{toda-b2}) satisfy (\ref{asm-u}), and  (\ref{ah}) holds for $h^k=(h_1^k,h_2^k)$. Then $(\sigma_1,\sigma_2)$ defined by (\ref{12sep2e3}) is one of the following types:
$$(2,0),(0,2),(4,2),(2,6),(4,8),(6,6),(6,8). $$
\end{thm}

As an application of Theoem \ref{cla-b2} we consider the $B_2$ Toda system defined on a compact Riemann surface $(M,g)$, whose volume is assumed to be $1 $ for convenience.
\begin{equation}\label{b2-m}
\left\{\begin{array}{ll}
\Delta_g u_1+2\rho_1(\frac{h_1e^{u_1}}{\int_Mh_1e^{u_1}dV_g}-1)-\rho_2(\frac{h_2e^{u_2}}{\int_Mh_2e^{u_2}dV_g}-1)=0,\\
\\
\Delta_g u_2-2\rho_1(\frac{h_1e^{u_1}}{\int_Mh_1e^{u_1}dV_g}-1)+2\rho_2(\frac{h_2e^{u_2}}{\int_Mh_2e^{u_2}dV_g}-1)=0.
\end{array}
\right.
\end{equation}
where $h_1,h_2$ are positive smooth functions on $M$, $\Delta_g$ is the Laplace-Beltrami operator $(-\Delta_g\ge 0)$, $\rho_1,\rho_2$ are positive constants.

It is well known that for solutions of the following general system on $M$:
$$\Delta_g u_i+\sum_{j=1}^n a_{ij}(\rho_j \frac{h_j e^{u_j}}{\int_M h_j e^{u_j} dV_g}-1)=0, \quad i=1,..,n, $$
the corresponding variational form is
$$\phi_{\rho}(u)=\frac 12 \sum_{i,j=1}^na^{ij}\int_M \nabla_g u_i \cdot \nabla_g u_j dV_g-\sum_{j=1}^n \rho_j \log \int_M h_j e^{u_j}dV_g. $$
where $(a^{ij})_{n\times n}=(a_{ij})_{n\times n}^{-1}$.

For solutions of (\ref{b2-m}) in \, \, $\hdot^{1,2}(M)$ we prove the following a priori estimate:
\begin{thm}\label{a-pri}
Suppose $\rho_1,\rho_2>0 $ and none of them is equal to a multiple of $4\pi$, then
$$ |u_i(x)|\le C, \quad i=1,2 $$
where $u=(u_1,u_2)$ is a solution of (\ref{b2-m}) in \, \, $\hdot^{1,2}(M)$.
\end{thm}

Next we consider the locally defined $G_2$ Toda system:

\begin{equation}\label{g2-loc}
\left\{\begin{array}{ll}
\Delta u_1^k+2h_1^ke^{u_1^k}-h_2^ke^{u_2^k}=0,\\
\\
\Delta u_2^k-3h_1^ke^{u_1^k}+2h_2^ke^{u_2^k}=0, \quad \mbox{ in } \quad B(0,1)
\end{array}
\right.
\end{equation}
where the assumptions on $u^k=(u_1^k,u_2^k)$ and $h^k=(h_1^k,h_2^k)$ are the same as in the $B_2$ case.  Then in this case we have

\begin{thm}\label{g2-loc-e}
Let $(\sigma_1,\sigma_2)$ be defined for blowup solutions $u^k$ of (\ref{g2-loc}) as (\ref{12sep2e3}). Suppose (\ref{asm-u}) holds for $u^k$, then $(\sigma_1,\sigma_2)$ satisfies
$$3\sigma_1^2+\sigma_2^2-3\sigma_1\sigma_2=6\sigma_1+2\sigma_2. $$
If in addition we have
$$\lim_{k\to \infty}\frac 1{2\pi}\int_{B(0,1)} h_1^k e^{u_1^k}<4+2\sqrt{2},\quad \lim_{k\to \infty}\frac 1{2\pi}\int_{B(0,1)} h_2^k e^{u_2^k}<10+2\sqrt{7}, $$
$(\sigma_1,\sigma_2)$ is one of the following types:
$$(2,0),(0,2),(4,2),(2,8),(4,12). $$
\end{thm}

As an application of Theorem \ref{g2-loc-e} we consider the following $G_2$ Toda system defined on the compact Riemann surface $(M,g)$:

\begin{equation}\label{g2-m}
\left\{\begin{array}{ll}
\Delta_g u_1+2\rho_1(\frac{h_1e^{u_1}}{\int_Mh_1e^{u_1}dV_g}-1)-\rho_2(\frac{h_2e^{u_2}}{\int_Mh_2e^{u_2}dV_g}-1)=0,\\
\\
\Delta_g u_2-3\rho_1(\frac{h_1e^{u_1}}{\int_Mh_1e^{u_1}dV_g}-1)+2\rho_2(\frac{h_2e^{u_2}}{\int_Mh_2e^{u_2}dV_g}-1)=0.
\end{array}
\right.
\end{equation}

Corresponding to Theorem \ref{a-pri} we have

\begin{thm}\label{a-pri-g2} Let $u=(u_1,u_2)$ be a solution of (\ref{g2-m}) in \,\, $\hdot^{1,2}(M)$, if $\rho_1\in (0,(8+4\sqrt{2})\pi),\rho_2\in (0, (20+4\sqrt{7})\pi) $ and none of them is equal to a multiple of $4\pi$, there is a constant $C$ independent of $u$ such that
$$ |u_i|\le C, \quad i=1,2. $$
\end{thm}

Theorem \ref{a-pri} and Theorem \ref{a-pri-g2} make it possible to calculate the Leray-Schauder degree $d_{\rho_1,\rho_2}$ for $(\rho_1,\rho_2)$ in intervals determined by multiples of $4\pi$ and
the unreasonable energy upper bound in (\ref{g2-m}). In order not to make this article exceedingly long we shall address the degree counting formula in a subsequent paper.

To end the introduction we describe the outline of the proof of Theorem \ref{cla-b2} and Theorem \ref{g2-loc-e}. First we use a selection process to determine a finite number of mutually disjoint bubbling disks. The idea of the selection process was first introduced by Schoen \cite{schoen} and is very useful for prescribing curvature type equations. In \cite{lwz-apde} Lin-Wei-Zhang applied this method to locate bubbling disks for systems of equations defined in two dimensional spaces. In each of the aforementioned bubbling disks, a partial blowup phenomenon occurs, which means if the blowup solutions are scaled according to the maximum of both components in the disk, only one component converges to a single equation after taking the limit. The phenomenon of partial blowup is the major difficulty for understanding the profile of bubbling solutions for systems. After identifying bubbling disks we use a Harnack inequality (Proposition B below) proved in \cite{lwz-apde} to describe the behavior of each component according to its spherical average around each blowup point. Each component is called to have ``fast decay" or ``slow decay" (see the next section for definition) based on its behavior. Roughly speaking, among the two components, at least one of them has fast decay and the energy of which is determined. The energy of the other component can be determined by the Pohozaev identity. Since we have to use the Pohozaev identity to determine the energy of one component, this approach only works for systems of two equations. Once the energy of at least one component is determined in each bubbling disk, we can also do the same for bubbling disks in a ``group" (see \cite{lwz-apde}). A group of bubbling disks looks roughly like a single disk after scaling. One major new ingredient in this article is to rule out the situation that there are only two bubbling disks in one group. For this we shall use Eremenko's work on surfaces with conical singularities to calculate $\frac 1{2\pi}\int_{\mathbb R^2}e^u$
where $u$ is a solution of
$$\Delta u+2e^u=4\pi \gamma_1 \delta_{p_1}+4\pi \gamma_2 \delta_{p_2}, \quad \int_{\mathbb R^2} e^u<\infty  $$
where $p_1,p_2$ are two disjoint points in $\mathbb R^2$ and $\gamma_1,\gamma_2>-1$.

We were not able to prove that all the energy types for $G_2$ Toda systems are multiples of $2$ without assumption. The reason that some strange numbers appear in Theorem
\ref{g2-loc-e} is because we don't have good estimates of the energy of Liouville solution with more than two singular points.

The organization of this article is as follows. In section two we mention a few tools we shall use in the proof of the main theorems. Some majors tools are developed in \cite{lwz-apde} and a new tool is based on Eremenko's work. Then in sections three and four we prove the concentration theorems and the proof of a priori estimates can be found in section five.
\medskip

\emph{Acknowledgement} Part of the work was finished when the second author was visiting Taida Institute for Mathematical Sciences (TIMS) in December 2014. He would like to express his deep gratitude to TIMS for their warm hospitality and financial aid.

\section{Preliminary discussions}
In this section we list a number of tools we shall use in the proof of Theorem \ref{cla-b2} and Theorem \ref{g2-loc-e}. Many of them come from the previous work of the authors and J. Wei \cite{lwz-apde} and a result of A. Eremenko \cite{eremenko}.

In \cite{lwz-apde} Lin-Wei-Zhang studied the concentration of energy for blowup solutions to
\begin{equation}\label{toda-an}
\Delta u_i+\sum_j a_{ij} h_je^{u_j}=4\pi \gamma_i \delta_0, \quad i=1,..,n, \mbox{ in } \, B(0,1)
\end{equation}
where $A=A_n$ is the  Cartan matrix of order $n$, $h_1...,h_n$ are positive smooth functions on $\bar B(0,1)$, $\gamma_i>-1$ indicates the strength of the Dirac mass at $0$.
Here we recall Proposition 2.1 of \cite{lwz-apde} which was the result of a selection process ( see \cite{schoen}):

\emph{ Proposition A:
Let $u^k=(u_1^k,..,u_n^k)$ be a sequence of solutions to (\ref{toda-an}) with $\gamma_1=..=\gamma_n=0$. Suppose there is a uniform bound for $\int_{B(0,1)}h_i^k e^{u_i^k}$ and the oscillation of
$u_i^k$ on $\partial B(0,1)$. Then if $0$ is the only blowup point of $u^k=(u_1^k,...., u_n^k)$,
there exist finite sequences of points $\Sigma_k:=\{x_1^k,....,x_m^k\}$ (all $x_j^k\to 0, j=1,...,m$) and positive numbers $l_1^k,...,l_m^k\to 0$ such that the following four properties hold:
\begin{enumerate}\item
$\max_{i\in I}\{ u_i^k(x_j^k)\}=\max_{B(x_j^k,l_j^k),i\in I}\{ u_i^k\}$ for all $j=1,..,m$, where $I=\{1,..,n\}$.
\item
$exp (\frac 12 \max_{i\in I}\{u_i^k(x_j^k)\}) l_j^k\to \infty, \quad j=1,...,m$.
\item There exists $C_1>0$ independent of $k$ such that
$$ u_i^k(x)+2\log \,\, dist(x,\Sigma_k)\le C_1, \quad \forall x\in B(0,1),\quad i\in I=\{1,...,n\} $$
where $dist$ stands for distance.
\item In each $B(x_j^k,l_j^k)$ let
\begin{equation}\label{vik-sel}
v_i^k(y)=u_i^k(\epsilon_k y+x_j^k)+2\log \epsilon_k, \quad \epsilon_k=e^{-\frac 12 M_k},\quad M_k=\max_i\max_{B(x_j^k,l_j^k)}u_i^k.
\end{equation}
 Then one of the following two alternatives holds\\
 (a):\quad The sequence is fully bubbling: along a subsequence $(v_1^k,...,v_n^k)$ converges in $C^2_{loc}(\mathbb R^2)$ to $(v_1,...,v_n)$ which satisfies
$$\Delta v_i+\sum_{j\in I} a_{ij}h_j e^{v_j}=0,\quad \mathbb R^2, \quad i\in I. $$
$$\lim_{k\to \infty} \int_{B(x_j^k,l_j^k)} \sum_{t\in I} a_{it}h_t^ke^{u_t^k}>4\pi,\quad i\in I.
$$
  (b):$I=J_1\cup J_2\cup... \cup J_m\cup N$ where $J_1,J_2,...,J_m$ and $N$ are disjoint sets. $N\neq \emptyset$ and each $J_t$ ($t=1,..,m$) consists of consecutive indices if it has more than one index.  For each $i\in N$,  $v_j^k$ tends to $-\infty$ over any fixed compact subset of $\mathbb R^2$. The components of $v^k=(v_1^k,...,v_n^k)$ corresponding to each $J_l$ ($l=1,...,m$) converge in $C^2_{loc}(\mathbb R^2)$ to a $SU(|J_l|+1)$ Toda system, where $|J_l|$ is the number of indices in $J_l$. For each $i\in J_l$, we have
  $$\lim_{k\to \infty}\int_{B(x_j^k,l_j^k)}\sum_{t\in J_l}a_{it}h_t^ke^{v_t^k}>4 \pi. $$
\end{enumerate}
}

\medskip

The selection process of Schoen singles out a finite number of bubbling disks for the sequence of blowup solutions $u^k=(u_1^k,...,u_n^k)$. In each of the bubbling disks, at least one component has energy greater than $2$:
$$\frac 1{2\pi}\int_{B(x_k,l_k)}h_i^k e^{u_i^k}>2.  $$
Here we use the integral above to denote the energy of $u_i^k$ in $B(x_k, l_k)$. Since there is a uniform bound of the energy of all components, there are only finite bubbling disks. The selection process also determines that
$$u_i^k(x)+2\log dist(x,\Sigma_k)\le C_1, \quad x\in B(0,1), $$
which provides a control on the upper bound of the behavior of blowup solutions outside the bubbling disks.

\medskip

The following proposition in \cite{lwz-apde} plays an essential role in the proof of main results in their article:

\emph{ Proposition B: (Lemma 2.1 of \cite{lwz-apde}) For all $x_0\in B(0,1)\setminus \Sigma_k$, there exists $C_0$ independent of $x_0$ and $k$ such that
$$|u_i^k(x_1)-u_i^k(x_2)|\le C_0,\quad \forall x_1,x_2\in B(x_0,d(x_0,\Sigma_k)/2),\mbox{ for all } i\in I. $$
}

Proposition B is a Harnack type estimate which reveals important information on the behavior of blowup solutions away from the bubbling area. Let $x_k\in \Sigma_k$ and $\tau_k=\frac 12 dist(x_k,\Sigma_k\setminus \{x_k\})$, then for $x,y \in B(x_k,\tau_k)$ and $|x-x_k|=|y-x_k|$ we have $u_i^k(x)=\bar u_{x_k,i}(r)+O(1)$ where
$r=|x_k-x|$ and
$$\bar u_{x_k,i}(r)=\frac 1{2\pi r}\int_{\partial B(x_k, r)}u_i^k dS. $$
In other words, the behavior of $u_i^k$ outside the bubbling disks can be represented by its spherical average in a neighborhood of a point in $\Sigma_k$.

For each $x_k\in \Sigma_k$, let
$$\sigma_{i,x_k}(r)=\frac 1{2\pi}\int_{B(x_k,r)}h_i^ke^{u_i^k}, \quad i=1..n, \quad r\le \tau_k=\frac 12 dist(x_k, \Sigma_k\setminus \{x_k\}), $$
we have
$$
\frac{d}{dr}\bar u_{x_k,i}(r)=\frac{1}{2\pi r}\int_{\partial B(x_k,r)}\frac{\partial u_i^k}{\partial \nu}=\frac{1}{2\pi r}\int_{B(x_k,r)}\Delta u_i^k=
-\sum_{j=1}^n\frac{a_{ij}\sigma_{j,x_k}(r)}r,
$$
For each component $i$ we say $u_i^k$ has fast decay on $x\in B(0,1)$ if
$$u_i^k(x)+2 \log dist(x,\Sigma_k)\le -N_k $$
for some $N_k\to \infty$. If there is a $C\in \mathbb R$ independent of $k$ and
$$u_i^k(x)+2 \log dist(x,\Sigma_k)\ge -C $$
we say $u_i^k$ has a slow decay at $x$. Here we note that
$$u_i^k(x)+2 \log dist(x,\Sigma_k)\le C $$
holds for all $x\in B(0,1)$.

The definition of fast and slow decay is very important for evaluating Pohozaev identities. For example in $B(x_k,r)$, the following Pohozaev identity holds for solutions to
$$\Delta u_i+\sum_{i=1}^n a_{ij}h_j^ke^{u_j^k}=0$$
with coefficient matrix being the Cartan matrix (in fact as long as $(a_{ij})$ is symmetric and invertible the following holds as well):
\begin{eqnarray}
&& \sum_i \int_{B(x_k,r)}(x\cdot \nabla  h_i^k)e^{u_i^k}+2\sum_i \int_{B(x_k,r)} h_i^k e^{u_i^k} \nonumber \\
&=& r\int_{\partial B(x_k,r)} \sum_i h_i^k e^{u_i^k}+r\int_{\partial B(x_k,r)}\sum_{i,j}\big ( a^{ij}\partial_{\nu}u_i^k\partial_{\nu} u_j^k-\frac 12 a^{ij}\nabla u_i^k \nabla u_j^k\big )
\label{poho1}
\end{eqnarray}
In order to evaluate the energy concentration from (\ref{poho1}) it is important to choose $r$ so that all components on $\partial B(x_k,r)$ have fast decay. Otherwise the first term on the right hand side of (\ref{poho1}) is not $o(1)$.

Finally we list a major new tool on the total energy of the following equation:
\begin{equation}\label{2-pole}
\left\{\begin{array}{ll}
\Delta u+2 e^u=4\pi(\theta_1-1)\delta_{p_1}+4\pi(\theta_2-1)\delta_{p_2}, \quad \mbox{ in }\quad \mathbb R^2, \\
\\
\int_{\mathbb R^2}e^u<\infty.
\end{array}
\right.
\end{equation}
Based on a theorem of Eremenko \cite{eremenko} we shall use the following theorem:

\emph{Theorem B-1: Let $p_1,p_2$ be two distinct points in $\mathbb R^2$ and $\theta_1,\theta_2$ be positive integers,
 then any solution $u$ of (\ref{2-pole}) satisfies
$$\frac 1{2\pi}\int_{\mathbb R^2}e^{u}=\theta_1+\theta_2+\theta_3-1 $$
where $\theta_3$ is a positive integer such that $\theta_1+\theta_2+\theta_3$ is odd and $\theta_i+\theta_j>\theta_k$ for $(i,j,k)$ being
any permutation of $(1,2,3)$.
}

\medskip

Theorem B-1 can be found in the appendix.

\section{Proof of Theorem \ref{cla-b2}}

First we observe that (\ref{toda-b2}) is a special case of $SU(4)$ Toda system by letting $u_3=u_1$ in $SU(4)$ Toda system. Consequently Proposition A can be applied to (\ref{toda-b2}) to obtain the
following blowup set $\Sigma_k=\{x_1^k,...,x_m^k\}$ and $l_1^k,...,l_m^k\in \mathbb R$ such that all $x_l^k$ ($l=1,...,m$) tend to the origin and all $l_i^k$ ($i=1,...,m$) tend to $0$. Moreover the following properties hold:

\begin{enumerate}
\item $\max_{i}\{u_i^k(x_j^k)\}=\max_{x\in B(x_j^k,l_j^k)}\{u_j^k(x)\}$, for all $j=1,...,m$.
\item $exp(\frac 12 \max_{i}\{u_i^k(x_j^k)\}l_j^k\to \infty,\quad j=1,...,m. $
\item There exists $C_1>0$ independent of $k$ such that
$$u_i^k(x)+2\log dist(x,\Sigma_k)\le C_1,\quad \forall x\in B(0,1), \quad i=1,2. $$
where $dist$ stands for distance.
\item In each $B(x_j^k,l_j^k)$ let $M_k=\max_i\max_{B(x_j^k,l_j^k)}u_i^k$, $\epsilon_k=e^{-\frac 12 M_k}$ and
$$v_i^k(y)=u_i^k(\epsilon_ky+x_j^k)+2\log \epsilon_k, \quad \mbox{ for } \quad \epsilon_ky+x_j^k\in B(x_j^k,l_j^k) $$
Then one of the following two alternatives holds:

(a): The sequence is fully bubbling: along a subsequence $(v_1^k,v_2^k)$ converges in $C^2_{loc}(\mathbb R^2)$ to $(v_1,v_2)$ that satisfies
$$\left\{\begin{array}{ll}
\Delta v_1+2e^{v_1}-e^{v_2}=0, \\
\Delta v_2-2e^{v_1}+2e^{v_2}=0,\quad \mbox{ in }\quad \mathbb R^2,\\
\int_{\mathbb R^2}e^{v_1}+e^{v_2}<\infty.
\end{array}
\right.
$$
$$\lim_{k\to \infty} \frac 1{2\pi}\int_{B(x_j^k,l_j^k)}\sum_{t=1}^2a_{it}h_t^ke^{u_t^k}>2, \quad i=1,2 $$
where $(a_{ij})$ is the $B_2$ Cartan matrix.

(b): $v^k=(v_1^k,v_2^k)$ has only one component converging to a single equation. If it is the first equation the limit equation
is
$$\Delta v_1+2e^{v_1}=0,\quad \mbox{ in }\quad \mathbb R^2. $$
If it is the second equation it is
$$\Delta v_2+2e^{v_2}=0,\quad \mbox{in}\quad \mathbb R^2. $$
In either case, the convergent component satisfies
$$\lim_{k\to \infty}\frac 1{2\pi}\int_{B(x_i^k,l_i^k)}2h_i^ke^{v_i^k}>2. $$
Here for convenience we assume $\lim_{k\to \infty} h_i^k(0)=1$. This assumption is not essential.
\end{enumerate}

\begin{rem} $l_i^k$ can be chosen so that for $x_j^k\in \Sigma_k$, let $t_k=dist(x_j^k,\Sigma_k\setminus \{x_j^k\})$, then $t_k/l_j^k\to \infty$.
\end{rem}

Since the $B_2$ Toda system is a special case of the $A_3$ Toda system, the Pohozaev identity for $A_3$ Toda system can also be applied for the $B_2$ Toda system. Here we recall that for $A^*=(a^*_{ij})_{3\times 3}$ being the $A_3$ Cartan matrix we have
$$\sum_{i,j=1}^3a^*_{ij}\sigma_i\sigma_j=4\sum_{i=1}^3\sigma_i. $$
Replacing $\sigma_3$ by $\sigma_1$ we have
\begin{equation}\label{ener-pi}
2\sigma_1^2-2\sigma_1\sigma_2+\sigma_2^2=4\sigma_1+2\sigma_2.
\end{equation}

Next we consider the energy of global solutions. By the classification theorem of Lin-Wei-Ye \cite{lin-wei-ye}, if $u=(u_1,u_2,u_3)$ is a global solution of $SU(3)$ Toda system with finite energy, then
$$\sum_{j=1}^3a^*_{ij}\int_{\mathbb R^2}e^{u_j}=8\pi,\quad i=1,2,3 $$
where $(a^*_{ij})$ is the $A_3$ Cartan matrix. Therefore in our case we let $u_1=u_3$ and then we have
$$\left\{\begin{array}{ll}
2\sigma_1-\sigma_2=4,\\
-2\sigma_1+2\sigma_2=4
\end{array}
\right.
\qquad \mbox{ which gives } \quad \sigma_1=6,\,\, \sigma_2=8. $$
It is also easy to verify that $(6,8)$ is the ``largest energy" because if $(6+t,8+s)$ also satisfies (\ref{ener-pi}) with $s,t\ge 0$. Then $s=t=0$.

To understand the concentration of energy we start from any fixed member of $\Sigma_k$. Say, $x_1^k\in \Sigma_k$ and let $\tau_{1,k}=distance(x_1^k,\Sigma_k\setminus \{x_1^k\})/2$.  Let
$$\sigma_{i,k}(r)=\frac 1{2\pi }\int_{B(x_1^k,r)}h_i^ke^{u_i^k},\quad i=1,2, \quad r\le \tau_{1,k} $$
and
$$ \bar u_i(r)=\frac{1}{2\pi r}\int_{\partial B(x_1^k,r)}u_i^k, \quad i=1,2. $$
Direct computation shows for $r\in (0,\tau_{1,k})$
\begin{align*}
\frac{d}{dr}\bar u_1(r)=\frac{-2\sigma_{1,k}(r)+\sigma_{2,k}(r)}r, \\
\frac{d}{dr}\bar u_2(r)=\frac{2\sigma_{1,k}(r)-2\sigma_{2,k}(r)}r.
\end{align*}
 The selection process guarantees that
$$u_i^k(x)+2\log |x-x_1^k|\le C, \quad |x-x_1^k|\le \tau_{1,k}, \quad i=1,2. $$
If both components have fast decay on $\partial B(x_1^k,r)$ ($r\in (0,\tau_{1,k})$), $(\sigma_{1,k}(r), \sigma_{2,k}(r))$ satisfies
\begin{equation}\label{pi-2}
2\sigma_{1,k}(r)^2-2\sigma_{1,k}(r)\sigma_{2,k}(r)+\sigma_{2,k}^2(r)=4\sigma_{1,k}(r)+2\sigma_{2,k}(r)+o(1).
\end{equation}
If we write (\ref{pi-2}) as
$$\sigma_{1,k}(r)(2\sigma_{1,k}(r)-\sigma_{2,k}(r)-4)+\sigma_{2,k}(r)(\sigma_{2,k}(r)-\sigma_{1,k}(r)-2)=o(1), $$
we see that for any $r$, if both components have fast decay, either $2\sigma_{1,k}(r)-\sigma_{2,k}\ge 4+o(1)$, which means
$$
\frac{d}{dr}(\bar u_1(r)+2\log r)\le -\frac{2+o(1)}r$$
or $\sigma_{2,k}(r)-\sigma_{1,k}\ge 2+o(1)$, which implies
$$\frac{d}{dr}(\bar u_2(r)+2\log r)\le -\frac{2+o(1)}r. $$
If one component (say $u_1^k$) satisfies
$$\frac{d}{dr}(\bar u_1(r)+2\log r)>0, $$
there is a possibility that for some larger radius $s$, $u_1^k$ becomes a slow decay component on $\partial B(x_1^k,s)$.

Next we consider the possible energy concentration types in $B(x_1^k,\tau_{1,k})$.
 Let $B(x_1^k,l_1^k)$ be a bubbling disk. Suppose in this disk, the first component converges according to the scaling of the maximum. If the region is close enough to the center of blowup ( the region still tends to infinity after scaling according to the maximum), we have
$$\sigma_{1,k}(l_k)=2+o(1),\quad \sigma_{2,k}(l_k)=o(1). $$
Now we consider the energy change from $B(x_1^k,l_k)$ to $B(x_1^k,\tau_{1,k})$. First we notice that on $\partial B(x_1^k,l_k)$, $\frac{d}{dr}(\bar u_2+2\log r)>0$, which means $u_2^k$ may become a
slow decay component when $r$ increases. So the first possibility in $B(x_1^k,\tau_{1,k})$ is that
$$\sigma_{1,k}(\tau_{1,k})=2+o(1),\quad \sigma_{2,k}(\tau_{1,k})=o(1), $$
which means $u_2^k$ does not change to be a slow decay component. It is proved in \cite{lwz-apde} that if no component changes to a slow decay component, the energy of each component only changes by $o(1)$ ( see Lemma 5.1 of \cite{lwz-apde}. Even though that lemma addresses $SU(3)$ Toda system but a very similar proof also applies to this case).

Since $\frac{d}{dr}(\bar u_2+2\log r)>0$, $u_2^k$ could become a slow decay component before $r$ reaches $\tau_{1,k}$. Suppose at some $s>r$,
$$\bar u_2^k(s)+2\log s\ge -C $$
for some $C>0$ very large. Here we observe that at this moment $u_2^k$ starts to increase its energy but the energy of $u_1^k$ barely changes because $\frac{d}{dr}(\bar u_1^k+2\log r)$ is still negative. If $\tau_{1,k}/s\to \infty$, which means $\tau_{1,k}$ is very large comparing with $s$, there is $N$ such that at $\partial B(x_1^k,Ns)$
\begin{align*}
\sigma_{2,k}(Ns)\ge 5, \quad \sigma_{1,k}(Ns)=2+o(1) \\
\bar u_1^k(Ns)+2\log (Ns)\le -N_k, \quad \mbox{ for some } N_k\to \infty\\
\frac{d}{dr}(\bar u_2+2\log r)|_{r=Ns}<0, \quad \frac{d}{dr}(\bar u_1+2\log r)|_{r=Ns}>0.
\end{align*}
In other words, from $r=s$ to $r=Ns$, the energy of $u_2^k$ increases and as a result, the derivative of $\bar u_2(r)+2\log r$ changes from positive to negative. But because of the Harnack inequality (Proposition B) $u_1^k$ is still a fast decay component and its energy barely changed, even though at $r=Ns$ the derivative of $\bar u_1+2\log r$ has become positive due to the change of the energy of $u_2^k$.
Since $\tau_{1,k}/s\to \infty$ we can find $N_k'$ tending to $0$ slowly such that $N_k's\le \tau_{1,k}/2$ and on $\partial B(x_1^k,N_k's)$ both $u_1^k$ and $u_2^k$ have fast decay. Evaluating
the Pohozaev identity on $\partial B(x_1^k,N_k's)$ we have
$$\sigma_{1,k}(N_k's)=2+o(1),\quad \sigma_{2,k}(N_k's)=6+o(1). $$
If $\tau_{1,k}$ is only comparable to $s$, then on $\partial B(x_1^k,\tau_{1,k})$, $u_2^k$ is a slow decay component and $\sigma_{1,k}(\tau_{1,k})=2+o(1)$.

At $\partial B(x_1^k,N_k's)$,
$$\frac{d}{dr}(\bar u_1^k(r)+2\log r)=\frac{4+o(1)}r, \quad r=N_k's, $$
and
$$\frac{d}{dr}(\bar u_2^k(r)+2\log r)=\frac{-6+o(1)}{r},\quad r=N_k's. $$
So at this radius there is a possibility that $u_1^k$ may become a slow decay component for larger $r$.

By exactly the same reason as before it is possible that $u_1^k$ increases to a slow decay component on $\partial B(x_1^k,\tau_{1,k})$, for which the second component has the energy
$\sigma_{2,k}(\tau_{1,k})=6+o(1)$; or $u_1^k$ finishes its transition of energy before $r$ reaching $\tau_{1,k}$: $\exists s_k\le \tau_{1,k}$ such that both components have fast decay on $s_k$ and
$$\sigma_{1,k}(s_k)=6+o(1),\quad \sigma_{2,k}(s_k)=6+o(1). $$
Similarly if $s_k$ is small compared to $\tau_{1,k}$ we could also have $(6,8)$ as the energy type in $B(x_1^k,\tau_{1,k})$. Since $(6,8)$ is the type of the energy a global solution has, there is no extra energy outside (see Theorem 4.1 and Theorem 4.2 of \cite{lwz-apde}).

If we start with the type $(0,2)$, which means in $B(x_1^k,l_k)$, $\sigma_{1,k}(l_k)=o(1)$ and $\sigma_{2,k}(l_k)=2+o(1)$. Then the following type may occur:
\begin{equation}\label{2-four}
(0,2),(4,2),(4,8),(6,8).
\end{equation}
More specifically if both components have fast decay on $\partial B(x_1^k,\tau_{1,k})$ then
$$(\sigma_{1,k}(\tau_{1,k}),\sigma_{2,k}(\tau_{1,k}))=(a+o(1),b+o(1)) $$
where $(a,b)$ is one of the four types in
(\ref{2-four}).

If one of the two components has slow decay on $\partial B(x_1^k,\tau_{1,k})$, from the discussion above we see that the energy of the other component (which has fast decay)
 is a multiple of $2$. For example, if the second component has slow decay on
$\partial B(x_1^k,\tau_{1,k})$, $\sigma_{1,k}(\tau_{1,k})$ is $2+o(1)$, $4+o(1)$ or $6+o(1)$. If the first component has slow decay, $\sigma_{2,k}(\tau_{1,k})$ is $2+o(1)$, $6+o(1)$ or $8+o(1)$.

Now we consider bubbling disks in a group. The concept of group is introduced in \cite{lwz-apde}), which means we consider bubbling disks relatively close to one another but relatively far away from other members in $\Sigma_k$. For example, if $x_1^k,x_2^k,x_3^k$ are in one group, $dist(x_i^k,x_j^k)\sim dist(x_m^k,x_n^k)$
 for $i,j,m,n\in \{1,2,3\}$ but $i\neq j$ and $m\neq n$ (here ``$\sim$" means comparable). Moreover for any $x_a^k\in \Sigma_k$ but $a\not \in \{1,2,3\}$, $dist(x_a^k,x_1^k)/dist(x_2^k,x_1^k)\to \infty$.
 By Proposition B if both components have fast decay around one bubbling disk, both components have fast decay around any of the disks in this group. Suppose the group members are $B(x_1^k,\tau_{1,k}),...,B(x_m^k,\tau_{m,k})$. By the definition of group, all the $\tau_{l,k}$ are comparable. If both components have
fast decay we can find $N_k\to \infty$ slowly such that all members in this group are contained in $B(x_1^k,N_k\tau_{1,k})$ and
\begin{align*}
\sigma_{1,k}(x_1^k,N_k\tau_{1,k})=\sum_{j=1}^m \sigma_{1,k}(x_j^k,\tau_{j,k})+o(1), \\
\sigma_{2,k}(x_1^k,N_k\tau_{1,k})=\sum_{j=1}^m \sigma_{2,k}(x_j^k,\tau_{j,k})+o(1),
\end{align*}
where
$$\sigma_{1,k}(x_1^k,N_k\tau_{1,k})=\frac 1{2\pi}\int_{B(x_1^k,N_k\tau_{1,k})}h_1^ke^{u_1^k} $$
and other notations are understood similarly.

In other words, $N_k$ can be chosen in a way that both components still have fast decay on $\partial B(x_1^k,N_k \tau_{1,k})$ and the energy of each component in this larger region is a $o(1)$ perturbation of the sum of the energy in each bubbling disk of this group. As observed before around each bubbling disk, at least one component has fast decay and the corresponding energy is a multiple of $2+o(1)$. Therefore both $\sigma_{1,k}(x_1^k,N_k\tau_{1,k})$ and $\sigma_{2,k}(x_1^k,N_k\tau_{1,k})$ are multiples of $2$ and they must satisfy the Pohozaev identity. The following are the only
pairs that satisfy the Pohozaev identity with each component being a multiple of $2$:
$$(2,0),(0,2), (4,2),(2,6),(4,8),(6,6),(6,8). $$

The final case we consider is when only one component has fast decay in a group. Suppose $u_1^k$ has fast decay in the group described as before. In this case
$\sum_{j=1}^m \sigma_{1,k}(x_j^k,\tau_{j,k}) $
is a multiple of $2$. In each bubbling disk, say $B(x_j^k,\tau_{j,k})$, since $u_2^k$ has slow decay, $\sigma_{2,k}(x_j^k,\tau_{j,k})>0$ and
by Proposition B we can choose $N_k\to \infty$ so that $u_1^k$ is still a fast decaying component on $\partial B(x_1^k, N_k \tau_{1,k})$ and
$$\sigma_{1,k}(x_1^k,N_k\tau_{1,k})=\sum_{j=1}^m \sigma_{1,k}(x_j^k,\tau_{j,k})+o(1). $$
Moreover,
$u_2^k$ also has
fast decay on $\partial B(x_1^k, N_k \tau_{1,k})$.  Therefore the Pohozaev identity can be evaluated on this radius.
Obviously this group is contained in $B(x_1^k, N_k\tau_{1,k})$.
Since the first component is a multiple of $2$, we see immediately that there is no new type except those we have known. For example, if we have two regions that both grow up from the type $(2,0)$, then the energy of the second component on $B(x_1^k,N_k \tau_{1,k})$ has to be $8+o(1)$.  If the group has two regions that grow from $(2,0)$ and $(4,2)$, the energy of the second component may grow to $8$ and make the energy type $(6,8)$. If there are types of $(2,0)$ in this group, the energy of $u_2^k$ on $B(x_1^k, N_k \tau_{1,k})$ has to be $6+o(1)$ to $8+o(1)$. If the energy of $u_1^k$ in this group is greater than $6$, for example $8+o(1)$, we use scaling to make the distance between any two blowup points in this group comparable to $1$. Then $u_2^k$, since it has slow decay in the neighborhood of these points, converges to a function $u_2$ that satisfies
\begin{equation}\label{limit-u2}
\Delta u_2+ 2 e^{u_2}=4\pi (\sum_{j=1}^m \gamma_{j}\delta_{p_j}), \quad \mbox{ in }\quad \mathbb R^2
\end{equation}
where $p_j$ ($j=1,..,m$) are the limits of blowup points in the group after scaling. Each $\gamma_j$ is a multiple of $2$ and $\sum_j\gamma_j=8$. By standard potential
analysis it is easy to see that there exist $\alpha>2$, $C\in \mathbb R$ and $\sigma>0$ such that
$$u_2(x)=-\alpha \log|x| +C+o(|x|^{-\sigma}), \quad |x|>1 $$
 and
$$\nabla u_2(x)=-\alpha \frac{x}{|x|^2}+O(|x|^{-\sigma-1}). $$
Integrating both sides of (\ref{limit-u2}) we see that
\begin{equation}\label{ener-u2-g}
\frac 1{2\pi}\int_{\mathbb R^2} e^{u_2}> 9.
\end{equation}
This estimate means in the group contained in $B(x_1^k, N_k\tau_{1,k})$ if $\sigma_{1,k}(x_1^k, N_k \tau_{1,k})=8+o(1)$,
$$\sigma_{2,k}(x_1^k, N_k \tau_{1,k})>9. $$
It is easy to see from direct computation that this pair of numbers cannot satisfy the Pohozaev identity.

If $u_2^k$ has fast decay in the group mentioned before, and $u_1^k$ has a slow decay in the neighborhood of the aforementioned bubbling disks, we first remark that if the total energy of $u_2^k$ in these bubbling disks is $6+o(1)$, the energy of $u_1^k$ in, say, $B(x_1^k, N_k \tau_{1,k})$ is either $2+o(1)$ or $6+o(1)$, according to
the Pohozaev identity. If the total energy of $u_2^k$ in this group is $8+o(1)$, the total energy of $u_1^k$ in $B(x_1^k, N_k \tau_{1,k})$ is either $4+o(1)$ or
$6+o(1)$. If the total energy of $u_2^k$ is $2m+o(1)$ for $m\ge 5, m\in \mathbb N$ (the set of natural numbers), as in the previous case we first derive a lower bound of the energy of $u_1^k$ in $B(x_1^k, N_k \tau_{1,k})$: By scaling the distance between any two members of the group is comparable to $1$, then the slow decaying component $u_1^k$ converges to $u_1$, which satisfies
\begin{equation}\label{limit-u1}
\Delta u_1+2 e^{u_1}=2\pi \sum_j \gamma_j \delta_{p_j}, \quad \mbox{ in }\quad \mathbb R^2.
\end{equation}
where $p_j$ are the images of blowup points in $\Sigma_k$ after scaling, $\sum_j \gamma_j=2m$ for $m\ge 5$. By standard potential analysis
$$\nabla u_1(x)=-\alpha \frac{x}{|x|^2}+O(|x|^{-\sigma-1}),\quad |x|>1 $$
for some $\alpha>2$ and $\sigma>0$. Integrating both sides of (\ref{limit-u1}) we have
$$\frac{1}{2\pi}\int_{\mathbb R^2} e^{u_1}>m+1. $$
Direct computation shows that there exists no pair $(\sigma_1,\sigma_2)$ satisfying the Pohozaev identity with $\sigma_2=2m$ ($m\ge 5$) and $\sigma_1>m+1$.

Finally we rule out the case that there are two $(0,2)$ type bubbling disks in the group and they are the only members. Note that this is the only case that the energy
of $u_2^k$ can be $4+o(1)$.
Since $u_1^k$ is slow decay. We see that by scaling the distance between these two groups into $1$, we see that the re-scaled version of $u_1^k$, which we use $u$ to represent, satisfies
$$\Delta u+2e^u=4\pi \delta_{p_1}+4\pi \delta_{p_2} \quad \mbox{ in } \quad \mathbb R^2. $$
and we also know $\int_{\mathbb R^2}e^{u}\le \infty$. By Theorem B-1 we know the total integration of
$$\frac 1{2\pi}\int_{\mathbb R^2}e^u=4 \mbox{ or  } 6. $$
But neither $(6,4)$ nor $(4,4)$ satisfies the Pohozaev identity, which means it is not possible to have two $(0,2)$ type regions in one group.

Since the combination of groups is similar to those of bubbling disks in one group, we have exhausted all the concentration types.
Theorem \ref{cla-b2} is established. $\Box$

\section{Partial results for $G_2$ Toda system}

The $G_2$ Toda system we consider is
\begin{equation}\label{e-g2}
\left\{\begin{array}{ll}
\Delta u_1+2h_1e^{u_1}-h_2 e^{u_2}=0, \\
\Delta u_1-3h_1 e^{u_1}+2h_2 e^{u_2}=0,\quad \mbox{in }\quad B(0,1)\subset \mathbb R^2
\end{array}
\right.
\end{equation}
and (\ref{ah}) and (\ref{asm-u}) also hold for $u^k=(u_1^k,u_2^k)$ as a sequence of solutions to (\ref{e-g2}).

The $G_2$ Toda system is a special case of a $A_6$ Toda system where $u_3=u_4=u_1+\log 2$, $u_5=u_2$ and $u_6=u_1$.  Using the Pohozaev identity for $A_6$ one obtains easily
the following identity:
\begin{equation}\label{pi-g2}
3\sigma_1^2+\sigma_2^2-3\sigma_1\sigma_2=6\sigma_1+2\sigma_2.
\end{equation}

First we consider global solution. We apply the theorem of Lin-Wei-Ye \cite{lin-wei-ye} to have
$$\sum_{j=1}^6a_{ij}\int_{\mathbb R^2}e^{u_j}=8\pi,\quad i=1,..,6 $$
where $(a_{ij})_{6\times 6}=A_6$. Using $\sigma_3=\sigma_4=2\sigma_1$, $\sigma_5=\sigma_2$, $\sigma_6=\sigma_1$, we see the blowup type for global solution is
$(12,20)$.

Next we consider partial blowup cases. Starting from $B(x_1^k, \tau_{1,k})$ which comes from the selection process as before. First we assume
 $$\sigma_{1,k}(l_k)=2+o(1), \quad \sigma_{2,k}(l_k)=o(1). $$
 From here we consider all the possible concentration types as $r$ increases from $l_k$ to $\tau_{1,k}$. If $u_2^k$ becomes a slow decaying component the energy type could change to $(2,8)$. After this $(8,8)$ could occur and finally $(8, 18)$ could occur.  If we start from $(0,2)$, then $u_1^k$ may change to a slow decaying component that leads to $(4,2)$. From here we could have $(4, 12)$ and $(10,12)$.

So far the following types are possibilities:

\begin{align}\label{all-g2}
(2,0),(0,2),(2,8),(4,2),(12,20),\\
(4,12),(8,8),(8,18),(10,12). \nonumber
\end{align}

When we consider the derivative of spherical averages, we have
$$ \frac d{dr}\bar u_1(r)=\frac{-2\sigma_{1,k}(r)+\sigma_{2,k}(r)}{r},\quad \frac{d}{dr} \bar u_2(r)=\frac{3\sigma_{1,k}(r)-2\sigma_{2,k}(r)}r. $$
If Pohozaev identity can be evaluated on $B(x, r)$, which means both $u_1^k$ and $u_2^k$ are fast decaying on $\partial B(x,r)$, either we have
\begin{equation}\label{g2-p1}
\bar u_1(r)\le -\frac{4+o(1)}r,\quad \mbox{ and } \quad \bar u_2(r)\ge -\frac{4+o(1)}r
\end{equation}
or
\begin{equation}\label{g2-p2}
\bar u_1(r)\ge -\frac{4+o(1)}r,\quad \mbox{ and } \quad \bar u_2(r)\le -\frac{4+o(1)}r
\end{equation}

On the other hand if both

\begin{equation}\label{g2-p3}
2\sigma_{1,k}(r)-\sigma_{2,k}(r)>2,\quad \mbox{ and }\quad -3\sigma_{1,k}(r)+2\sigma_{2,k}(r)>2
\end{equation}
on a fast decaying radius for both components, then there is no essential energy outside.

Next we consider the combination of bubbling disks in a group. If both components have fast decay in the neighborhood of those bubbling disks, both
$\sigma_{1,k}(x_1^k, N_k \tau_{1,k})$ and $\sigma_{2,k}(x_1^k, N_k \tau_{1,k})$ are multiples of $2$ and all such pairs that satisfy the Pohozaev identity have been listed in (\ref{all-g2}).  So we only consider the case that one component has slow decay in the neighborhood of those bubbling disks in one group. If the group has only two disks, we shall use Theorem B-1. From Theorem B-1 we see that the energy of the slow decaying component is also a multiple of $2$. So we still do not add any new type except those listed in (\ref{all-g2}).

Finally we rule out the case that there are at least three bubbling disks in one group. The most ``energy efficiency"  case is there are three $(2,0)$ type bubbling disks in one group or three $(0,2)$ type bubbling disks in one group. In the first case, clearly $u_2^k$ has slow decay and it is easy to see that
$$(6, 10+2\sqrt{7})$$
satisfies the Pohozaev identity. By assuming $\sigma_2<10+2\sqrt{7}$ we ruled out this case. If $u_1^k$ has the slow decay and there are three $(0,2)$ type bubbling disks in the group, we see easily that $(4+2\sqrt{2},6)$ satisfies the Pohozaev identity. By assume $\sigma_1<4+2\sqrt{2}$ this case is also ruled out. Other cases of having more bubbling disks or only three bubblings disks with more energy can also be ruled out easily by the restriction of $\sigma_i$.

Theorem \ref{g2-loc-e} is established. $\Box$

\section{ Proof of Theorem \ref{a-pri} and Theorem \ref{a-pri-g2}}

Since the nature of these two theorems is so close we just prove Theorem \ref{a-pri} as an example. The proof of Theorem \ref{a-pri-g2} is completely similar.

 Let
$$\tilde u_i=u_i-\log \int_M h_i e^{u_i},\quad i=1,2. $$
Then we have
\begin{equation}\label{b2-m1}
\left\{\begin{array}{ll}
\Delta_g \tilde u_1+2\rho_1(h_1e^{\tilde u_1}-1)-\rho_2 (h_2 e^{\tilde u_2}-1)=0,\\
\\
\Delta_g \tilde u_2-2\rho_1(h_1e^{\tilde u_1}-1)+2\rho_2 (h_2e^{\tilde u_2}-1)=0,
\end{array}
\right.
\end{equation}
and
\begin{equation}\label{ener-1}
\int_M h_i e^{\tilde u_i}dV_g=1,\quad i=1,2.
\end{equation}

In the first step, which is the major one we prove that there is a $C$ independent of $u$ such that
\begin{equation}\label{t-u-bound}
|\tilde u_1|+|\tilde u_2|\le C.
\end{equation}

In order to prove (\ref{t-u-bound}) we only need to prove
\begin{equation}\label{t-u-upp}
\max_i\max_M \tilde u_i\le C.
\end{equation}
because once (\ref{t-u-upp}) is established the lower bound can be obtained by standard Harnack inequality. Thus we shall prove
(\ref{t-u-upp}) by way of contradiction:
 Suppose $\tilde u^k=(\tilde u_1^k, \tilde u_2^k)$ is a sequence such that $\max_i\max_M\tilde u_i^k\to \infty$. Let $\mathbb G$ be the
blowup set for $\tilde u^k$:
$$\mathbb G=\{p\in M;\quad \exists x_k\to p, \quad \lim_{k\to \infty} \max_i \tilde u_i^k(x_k)\to \infty. \quad \} $$

Let $p_1$,...,$p_N$ be the blowup points of $\tilde u^k=(\tilde u_1^k,\tilde u_2^k)$ on $M$. It is easy to prove that there are only finite blowup points on $M$ by standard estimates ( see \cite{bm,lin-zhang-jfa,lin-wei-yang}, etc). Here we assume that $p_1,..,p_N$ are distinct points. By the Green's representation of $\tilde u_i^k$ it is easy to see that $\tilde u_i^k$ has bounded oscillation outside the bubbling area: $\cup_{i=1}^N B(p_i,\delta_0)$ where $\delta_0>0$ is chosen small enough to make
$$B(p_i,\delta_0)\cap B(p_j,\delta_0)=\emptyset.$$

Around each $p_i$, we write the equation in Euclidean form in local coordinates, let
$$\sigma_i(p_1):=\lim_{\delta\to 0}\lim_{k\to \infty} \frac 1{2\pi}\int_{B(p_1,\delta)}\rho_i^k h_i e^{\tilde u_i^k} dV_g,$$
then we shall see that this is exactly the same as what is defined for the local equation. Indeed, let $p_1^k\to p_1$ be where $\max_i\max_{B(p_1,\delta_0)}\tilde u_i^k$ is attained.
Take the local coordinates around $p_1^k$, then $ds^2$ has the form $e^{\psi(y_{p_1^k})}(dy_1^2+dy_2^2)$ where
$$|\nabla \psi(0)|=\psi(0)=0,\quad \Delta \psi=-2Ke^{\psi}, $$
$K$ is the Gauss curvature.
Obviously the equation for $\tilde u^k$ can be written as
\begin{equation}\label{t-u-2}
\left\{\begin{array}{ll}
\Delta \tilde u_1^k+2\rho_1^ke^{\psi}h_1e^{\tilde u_1^k}-\rho_2^ke^{\psi}h_2e^{\tilde u_2^k}+e^{\psi}(-2\rho_1^k+\rho_2^k)=0,\\
\\
\Delta \tilde u_2^k-2\rho_1^ke^{\psi}h_1e^{\tilde u_1^k}+2\rho_2^ke^{\psi}h_2e^{\tilde u_2^k}+e^{\psi}(2\rho_1^k-2\rho_2^k)=0, \quad \mbox{ in } \quad B_{\delta_0}
\end{array}
\right.
\end{equation}

Let $f_1^k$ and $f_2^k$ be defined by
$$ \left\{\begin{array}{ll}
\Delta f_1^k+e^{\psi}(-2\rho_1^k+\rho_2^k)=0,\quad \mbox{ in }\quad B_{\delta_0} ,\\
f_1^k=0,\quad \mbox{ on }\quad \partial B_{\delta_0}.
\end{array}
\right.
$$
and
$$ \left\{\begin{array}{ll}
\Delta f_2^k+e^{\psi}(2\rho_1^k-2\rho_2^k)=0,\quad \mbox{ in }\quad B_{\delta_0} ,\\
f_2^k=0,\quad \mbox{ on }\quad \partial B_{\delta_0}.
\end{array}
\right.
$$
By setting $\hat h_i^k=\rho_i^kh_ie^{\psi+f_i^k}$ and $\hat u_i^k=\tilde u_i^k+f_i^k$ we can write (\ref{t-u-2}) as
$$
\left\{\begin{array}{ll}
\Delta \hat u_1^k+2\hat h_1^k e^{\tilde u_1^k}-\hat h_2^k e^{\hat u_2^k}=0,\quad \mbox{in }\quad B_{\delta_0}, \\
\\
\Delta \hat u_2^k-2\hat h_1^k e^{\tilde u_1^k}+2\hat h_2^k e^{\hat u_2^k}=0,
\end{array}
\right.
$$
Since $dV_g=e^{\psi}dy$ it is easy to see that
$$\int_{B_{\delta_0}}\hat h_i^ke^{\hat u_i^k}dy=\int_{B(p_1^k,\delta_0)}\rho_i^k h_ie^{\tilde u_i^k}dV_g,\quad i=1,2. $$

Then $(\sigma_1(p_1),\sigma_2(p_1))$ is just one of the seven types. Moreover, at least one of the components has fast decay. For example if $(\sigma_1(p_1),\sigma_2(p_1))=(4,2)$, the second component decays fast. Since the oscillation of either component is finite away from the bubbling disks, we see that at least one component has little energy outside the bubbling disks. Also if the second component is fast decaying, the second component around any blowup point is also fast decaying. Then we have
$$\cup_{t=1}^N \frac 1{2\pi}\int_{B(p_t, \delta_0)} \rho_2^k h_2e^{\tilde u_2^k}dV_g=2 N $$
for some positive integer $N$. However since
$$\int_{M\setminus \cup B(p_t, \delta)}\rho_2^k h_2e^{\tilde u_2^k}dV_g=o(1) $$
and (\ref{ener-1}) holds, we get a contradiction to our assumption that $\rho_2^k$ cannot tend to a multiple of $4\pi$.
Thus (\ref{t-u-upp}), and consequently (\ref{t-u-bound})  are established.

Now we finish the proof of Theorem \ref{a-pri}. Clearly we can write the equation for $u=(u_1,u_2)$ as
$$\left\{\begin{array}{ll}
\Delta_g u_1+2\rho_1 (h_1 e^{\tilde u_1}-1)-\rho_2( h_2e^{\tilde u_2}-1)=0,\\
\\
\Delta_g u_2-2\rho_1 (h_1 e^{\tilde u_1}-1)+2\rho_2( h_2e^{\tilde u_2}-1)=0.
\end{array}
\right.
$$
Since
\begin{align*}
\int_M \bigg (2\rho_1 (h_1 e^{\tilde u_1}-1)-\rho_2( h_2e^{\tilde u_2}-1) \bigg )dV_g=0 ,\\
\int_M \bigg (-2\rho_1 (h_1 e^{\tilde u_1}-1)+2\rho_2( h_2e^{\tilde u_2}-1) \bigg )dV_g=0
\end{align*}
and (\ref{t-u-bound}) holds, from standard elliptic estimate and
$\int_M u_i=0$ ($i=1,2$) we see that $|u_1|+|u_2|\le C$.   Theorem \ref{a-pri} is established.
The proof of Theorem \ref{a-pri-g2} is similar. $\Box$

\section{Appendix: Two Theorems of Eremenko}

In this appendix we interpret some theorems of A. Eremenko \cite{eremenko} into a form that can be used in this article.

First we recall that on a Riemann surface $S$, a metric $g_0$ is called conformal if in any local coordinate system $z_l\in \Omega\subset \mathbb C$,
$$g_0(z_l)=e^u|dz_l|^2,\quad z_l\in \Omega $$
for a measurable and bounded function $u$ in $\Omega$. Let $P_0\in S$, a conformal metric $\tilde g_s$ is called to have a conical singularity at $P_0$ of total angle $2\pi(\alpha+1)$ ($\alpha>-1$) if there exist local coordinates $z(P)\in \Omega\subset \mathbb C$ and $u\in C^0(\Omega)\cap C^2(\Omega\setminus \{P\})$ such that $z(P_0)=0$ and
$$g_s(z)=|z|^{2\alpha}e^u |dz|^2, \quad z\in \Omega, $$
where $g_s$ is the local expression of $\tilde  g_z$.

In \cite{eremenko} A. Eremenko studied the following situations: Let  $p_1,p_2,p_3$ be distinct points on $\mathbb S^2$ and $2\pi \theta_j $($j=1,2,3$)  be their total angles of conic singularity respectively $(\theta_i>0,i=1,2,3)$.

\medskip
The first case is \emph{at least one of the three numbers $\theta_1,\theta_2,\theta_3$ is an integer}.  Among other things Eremenko proved:
\medskip

\emph{Theorem A (Eremenko): Let $\theta_1>1$ be an integer, but $\theta_2>0$ and $\theta_3>0$ are not integers. If there is a conformal metric with curvature $1$ on the sphere, with three conic singularities of angles $2\pi \theta_i$ ($i=1,2,3$), then either $\theta_2+\theta_3$ or $|\theta_2-\theta_3|$ is an integer $m$ of opposite parity from $\theta_1$, and $m\le \theta_1-1$. }

\medskip
\begin{rem} Two integers $A$ and $B$  are called to belong to opposite parity if one of them is even and the other is odd.
\end{rem}

Theorem A is the existence part of Eremenko's original statement in \cite{eremenko}, the uniqueness part of conformal metric with prescribed conical singularities can be found in  Fujimori, et. al \cite{fujimori}.

If all $\theta_i$ are positive integers not equal to $1$, Eremenko proved:

\medskip

\emph{ Theorem B (Eremenko): Let $\theta_1,\theta_2,\theta_3$ be three positive integers not equal to $1$. If there exists a conformal metric of curvature $1$ on the sphere, with three conic singularities of angles $2\pi \theta_1$, $2\pi \theta_2$ and $2\pi \theta_3$, respectively, then $\theta_+\theta_2+\theta_3$ is odd and
$\theta_i+\theta_j>\theta_k$ for $(i,j,k)$ being any permutation of $(1,2,3)$. }

\medskip

The second case is \emph{none of $\theta_1$,$\theta_2$,$\theta_3$ is an integer}. This is case is not used in this article but we still translate it into a PDE result for applications in the future.

\medskip

We say $(\theta_1,\theta_2,\theta_3)$ ( each $\theta_i>0$) is equivalent to $(\pm \theta_1+m,\pm \theta_2+n, \pm \theta_3+k)$ when $(m,n,k)$ are integers with the property $m+n+k=0$ ( mod $2$). Every non-integer triple is equivalent to one and only one triple with the property
\begin{equation}\label{uniq-tri}
0<\theta_1'+\theta_2'\le 1,\quad 0<\theta_2'+\theta_3'\le 1,\quad 0<\theta_1'+\theta_3'\le 1.
\end{equation}

For this case Eremenko proved:

\medskip

\emph{Theorem C (Eremenko): If none of $\theta_1,\theta_2,\theta_3$ is an integer, then a conformal metric of constant positive curvature on the sphere with conic singularities of total angles
$2\pi \theta_1$, $2\pi \theta_2$ and $2\pi \theta_3$ exists if and only if the unique equivalent triple with the property (\ref{uniq-tri}) satisfies $\theta_1'+\theta_2'+\theta_3'>1$. Such a metric of curvature $1$ is unique. }

\medskip

We shall interpret Theorems A,B and C for the following equation:
\begin{equation}\label{e-2s}
\Delta u+2e^u=4\pi(\theta_1-1)\delta_{p_1}+4\pi(\theta_2-1)\delta_{p_2},\quad \mbox{in}\quad \mathbb R^2,\quad \int_{\mathbb R^2}e^u<\infty.
\end{equation}
where $p_1,p_2$ are two distinct points in $\mathbb R^2$ and $\theta_1,\theta_2$ are positive constants.

\medskip

\medskip

\emph{Theorem A-1: Let $p_1,p_2$ be two distinct points in $\mathbb R^2$ and  $u$ be a solution of (\ref{e-2s}). Suppose $\theta_1$ is a positive integer and $\theta_2>0$ is not an integer, then any solution $u$  of (\ref{e-2s}) satisfies
$$\frac 1{2\pi}\int_{\mathbb R^2}e^u=\theta_1+\theta_2+\theta_3-1 $$
for some $\theta_3>0$. Moreover either $\theta_2+\theta_3$ or $|\theta_2-\theta_3|$ is an integer $m$ of opposite parity from $\theta_1$, and $m\le \theta_1-1$.
}

\medskip

The following theorem is a translation of Theorem B:

\medskip

\emph{ Theorem B-1: Let $p_1,p_2$ be two distinct points in $\mathbb R^2$ and $u$ be a solution of (\ref{e-2s}). If both $\theta_1$ and $\theta_2$ are positive integers. $u$ satisfies
$$\frac 1{2\pi} \int_{\mathbb R^2} e^u=\theta_1+\theta_2+\theta_3-1 $$
where $\theta_3$ is a positive integer, $\theta_1+\theta_2+\theta_3$ is odd and $\theta_i+\theta_j>\theta_k$ for $(i,j,k)$ being any permutation of $(1,2,3)$.
}

\medskip

The interpretation of Theorem C (which we don't use in this article)  is

\medskip

\emph{Theorem C-1: Suppose $\theta_1,\theta_2,\theta_3$ are all positive but none of them is an integer. Then there exists a unique solution $u$ of (\ref{e-2s}) with
$$\frac 1{2\pi}\int_{\mathbb R^2}e^u=\theta_1+\theta_2+\theta_3-1 $$
if and only if $\theta_1'+\theta_2'+\theta_3'>1$.
}

\bigskip

\noindent{\bf Proof of Theorem A-1,Theorem B-1 and Theorem C-1:}

\medskip

\noindent{\bf Step one}

In the first part of the proof, we state one fact: Let $q_1$, $q_2$ and $q_3$ be distinct points on $\mathbb C$ and
$2\pi\theta_i$, ($i=1,2,3$) be the total angle of conical singularity at $q_i$. Then the equation for
Gauss curvature equal to $1$ is
\begin{equation}\label{a-e1}
\Delta v+e^{2v}=2\pi(\theta_1-1)\delta_{q_1}+2\pi(\theta_2-1)\delta_{q_2}+2\pi(\theta_3-1)\delta_{q_3}.
\end{equation}
Here is the reason why (\ref{a-e1}) holds. First away from the singularities
$$K=-(\Delta v)e^{-2v}=1. $$
 Then around each singularity, say $q_1$, if in the neighborhood we let
$$\tilde v=v-(\theta_1-1)\log |x-q_1|, $$
the equation for $\tilde v$ around $q_1$ would be
$$\Delta \tilde v+e^{2\tilde v}|x-q_1|^{2(\theta_1-1)}=0. $$
Thus $2\pi\theta_1$ is the total angle at $q_1$.
Let $v_1=2v$, then (\ref{a-e1}) becomes
\begin{equation}\label{a-e2}
\Delta v_1+2e^{v_1}=4\pi(\theta_1-1)\delta_{q_1}+4\pi(\theta_2-1)\delta_{q_2}+4\pi(\theta_3-1)\delta_{q_3}.
\end{equation}

Theorem A and Theorem B can be applied if $\infty$ is a not a singular point. $\infty$ is not a singular point if
$$v_1(y)=-4\log |y|+O(1), \quad |y|>1 $$
because if we let $v_2(z)=v_1(z/|z|^2)-4\log |z|$ for $z$ close to the origin, we would have
$$\Delta v_2(z)+2e^{v_2(z)}=0  $$
in a neighborhood of $0$ and $v_2$ is bounded near $0$.

\medskip

\noindent{\bf Step two: we can assume $p_1=(-1,0)$, $p_2=(1,0)$:}

Going back to equation (\ref{e-2s}),
we first use a transformation to move $p_1$,$p_2$ in (\ref{e-2s}) to $(-1,0)$ and $(1,0)$: Let
$$u_1(x)=u(x)-2(\theta_1-1)\log |x-p_1|-2(\theta_2-1)\log |x-p_2|, $$
then clearly
$$\Delta u_1(x)=\Delta u-4\pi (\theta_1-1)\delta_{p_1}-4\pi(\theta_2-1)\delta_{p_2}, \quad \mbox{in  }\quad \mathbb R^2. $$
and
$$\Delta u_1+2e^{u_1}|x-p_1|^{2(\theta_1-1)}|x-p_2|^{2(\theta_2-1)}=0,\quad \mbox{in }\quad \mathbb R^2. $$
Let
$$u_2(y)=u_1(\frac{p_1+p_2}2+dy)+2(\theta_1+\theta_2-1)\log d, \quad d=|p_1-p_2|. $$
Then it is easy to verify
$$\Delta u_2(y)+2|y-P_1|^{2(\theta_1-1)}|y-P_2|^{2(\theta_2-1)}e^{u_2}=0,\quad \mbox{in }\quad \mathbb R^2
$$
where $P_1=(-1,0)$, $P_2=(1,0)$. Clearly $0$ is a regular point for $u_2$.

\medskip

\noindent{\bf Step Three: Singularity at infinity}

Let
$$u_3(y)=u_2(y)+2(\theta_1-1)\log |y-P_1|+2(\theta_2-1)\log |y-P_2|, $$
we have
$$\Delta u_3=\Delta u_2+4\pi(\theta_1-1)\delta_{P_1}+4\pi(\theta_2-1)\delta_{P_2}$$
and
\begin{equation}\label{e-u3}
\Delta u_3+2e^{u_3}=4\pi(\theta_1-1)\delta_{P_1}+4\pi(\theta_2-1)\delta_{P_2}.
\end{equation}
Then we consider the Kelvin transformation of $u_3$:
$$u_4(z)=u_3(\frac{z}{|z|^2})-4\log |z|, \quad z\in \mathbb R^2, $$
then round $0$,
$$\Delta u_4(z)+2e^{u_4}=0,\quad \mbox{ in }\quad B_{1/2}\setminus \{0\}. $$
In order to determine $\Delta u_4(0)$ we first observe that by standard potential analysis
\begin{equation}\label{a-e3}
\left\{\begin{array}{ll}
u_3(y)=-\alpha \log |y|+ c_1+O(|y|^{-\sigma}),\quad |y|>1,\\
\\
\nabla u_3(y)=-\alpha \frac{y}{|y|^2}+O(|y|^{-\sigma-1}), \quad |y|>1.
\end{array}
\right.
\end{equation}
for some $\alpha>2$, $c_1\in \mathbb R$ and $\sigma>0$.

By the definition of $u_3$ and $u_4$ we see that (\ref{a-e3}) leads to
$$u_4(z)=(\alpha -4)\log |z|+O(1) \quad \mbox{ near } \quad 0 $$
and
$$\Delta u_4(0)=2\pi (\alpha-4)\delta_0=4\pi(\frac{\alpha}2-2)\delta_0. $$
If we use $2\pi\theta_3$ to denote the total angle at $0$ we have
\begin{equation}\label{theta-3}
\theta_3=\frac{\alpha}2-1.
\end{equation}

\medskip

\noindent{\bf Step four: Completion of the proof:}

Now we integrate on both sides of (\ref{e-u3}), the left hand side gives (using (\ref{a-e3}))
$$\lim_{R\to \infty}\int_{B_R}\Delta u_3=-2\pi \alpha=-2\pi(2\theta_3+2) $$
where the last step is by (\ref{theta-3}). Direct computation shows that
$$\int_{\mathbb R^2}e^{u_3}=\int_{\mathbb R^2}e^u. $$
Thus
$$-2\pi(2\theta_3+2)+2\int_{\mathbb R^2}e^u=4\pi(\theta_1-1)+4\pi(\theta_2-1), $$
and we have obtained the following important equation:
\begin{equation}\label{major-e}
\frac 1{2\pi}\int_{\mathbb R^2}e^u=\theta_1+\theta_2+\theta_3-1.
\end{equation}
To finish the proof of Theorem A-1, we see that it follows directly from Theorem A if $\theta_1\neq 1$. If $\theta_1=1$, it is easy to verify that Theorem A-1 still holds, because the requirements on $\theta_2$ and $\theta_3$ imply $\theta_2=\theta_3$ and it is easy to check that the conclusion still holds by Prajapat and Tarentello's classification theorem \cite{prajapat}. See \cite{tr-2} as well. Theorem B-1 follows from Theorem B if no $\theta_i$ is $1$. If some $\theta_i$ is one, say $\theta_1=1$, the classical result of Troyanov \cite{tr-2} asserts $\theta_2=\theta_3$. It is easy to see that the conclusion of Theorem B-1 still holds in this case. Finally for this case we note that $\theta_1=\theta_2$ if $\theta_3=1$.
Theorem C-1 is straight forward interpretation of Theorem C.
 $\Box$

\end{document}